\documentclass[12pt,amscd]{amsart}
\usepackage{amssymb}  	

\usepackage{setspace}	

\input xypic \xyoption{curve}	

\usepackage{graphicx} 

\newcommand{\B}[1]{{\mathbf#1}} 
\newcommand{\C}[1]{{\mathcal#1}} 

\theoremstyle{plain}
\newtheorem{theorem}{Theorem}[section]

\newtheorem{prop}{Proposition}[section]

\theoremstyle{definition}
\newtheorem{defin}{Definition}[section]

\theoremstyle{definition}
\newtheorem{example}{Example}[section]

\theoremstyle{remark}
\newtheorem{rem}{Remark}[section]

\usepackage{fancyhdr}
\begin{document}

\title{Lattice Models of Finite Fields}
\author{Lucian M. Ionescu, Mina M. Zarrin} 
\address{Department of Mathematics, Illinois State University, IL 61790-4520}
\email{LMIones@ilstu.edu}
\email{MZarrin@ilstu.edu}
\date{July 2017} 

\begin{abstract}
Finite fields form an important chapter in abstract algebra, and mathematics in general. 
We aim to provide a geometric and intuitive model for finite fields, involving algebraic numbers,
in order to make them accessible and interesting to a much larger audience.

Such lattice models of finite fields provide a good basis for later developing the theory in a more concrete way,
including Frobenius elements, all the way to Artin reciprocity law.

Examples are provided, intended for an undergraduate audience in the first place.
\end{abstract}

\maketitle
\setcounter{tocdepth}{3} 
\tableofcontents

\section{Introduction}       
Finite fields are important mathematical structures, taking the learner from the familiar realm of congruence arithmetic
to algebraic number theory territory, and providing new tools for mathematical physics and cryptography,
for example.

We aim to highlight a pedagogical tool for the introduction of higher dimensional finite fields, 
which balances the traditional ``axiomatic'', top-down approach of Abstract Algebra, 
with a constructive, yet intuitive approach, using so called lattice models.
The ``standard'' way extends the primary finite field $F_p$, 
as a quotient of a polynomial ring $F_p[X]/<f(X)>$.
The lattice model approach extends the lattice of integers first, 
to place it in the context of complex numbers,
followed by the quotient modulo a prime.
In this way it mimics the elementary case of primary finite fields $F_p=Z/pZ$,
providing also a geometric intuition accompanied by the corresponding analytic-topologic tools
available.
$$
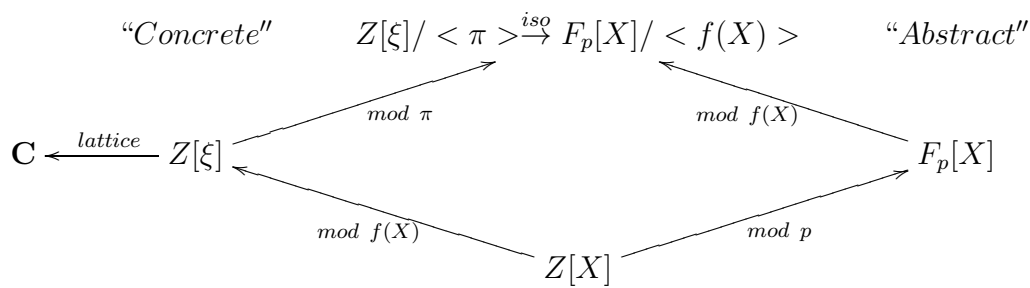
\diagram\label{Diagram1}
& ``Concrete'' & Z[\xi]/<\pi> \overset{iso}{\to}  F_p[X]/<f(X)>  & ``Abstract''\\
\B{C} & \lto_{lattice} Z[\xi]  \urto_{mod \ \pi} & & F_p[X] \ulto^{mod \ f(X)} & \\
&  & Z[X] \ulto^{mod\ f(X)} \urto_{mod \ p} & 
\enddiagram$$

Comparing with the concept of group, the ``abstract way'' is to define the algebraic structure with one binary operation, 
and then derive their properties from ``axioms'', perhaps too soon, 
before the student has enough examples to develop the ``feeling'' and intuition of what they are. 
The two dual, symbiotic types of groups, are the non-commutative groups of transformations, which always act on some space, 
and those we call Abelian, which in fact are ``discrete vector spaces'' on which the first kind act upon.
The ``unified'' approach through generalization and abstraction has its price: treating alike the two becomes the norm, 
and the differences in interpretation neglected. 

In this modern algebraic way of introducing algebraic structures abstractly, through general definitions,
and then quickly deriving their properties,
one would immediately ask the question of existence and uniqueness.
The later can be addressed in complete generality, without even knowing if they exist.
Existence is proven, of course, by constructing finite fields explicitly.

For primary finite fields $F_p$, in characteristic $p$, this is easy:  the well known Abelian groups $Z/nZ$, 
taught while doing congruence arithmetic, or rather viewed as rings $Z/nZ$,
are easily shown to be fields, when $n=p$ is a prime number;
but the other high dimensional finite fields are $F_p^d$ are harder to construct, 
and the ``future algebraist'', the student, learns by heart the recipes for constructing field extensions.

Pedagogically, examples should be provided first, 
worked with them to the point the student begins to like them, 
and then ``frame them'' in the appropriate axiomatic context.

The lattice models of finite field presented in this paper represent construction of $F_p^d$, 
generalizing the above simple case of wrapping the 1-dimensional lattice $Z$, with period 
corresponding to the prime ideal $pZ$.
By using higher dimensional lattice, instead of the standard adjunction of ``roots'' construction,
we provide a geometric interpretation, 
with a corresponding graphical representation which brings geometry up-front, 
to enjoy and play with ... if time allows it!

Of course, there is a price to pay: some new number systems need to be introduced along the way, 
still extensions using the standard algebraic construction, 
but so important that they need to be made well known well before the theory of finite fields takes off:
Gaussian and Eisenstein integers, and their generalizations (cyclotomic extensions).

And yet here again, one can borrow the geometric interpretation of complex numbers as representing
2D-rotations, and still provide enough geometric intuition, to overcome the abstract ``magical act'' 
of adjoining new symbols;
at least this is the opinion of one of the authors.

The article is organized as follows.
The next section \S\ref{S:BG}, introduces finite fields abstractedly, as in most textbooks of abstract algebra.
Section \S\ref{S:LM} constructs finite fields as congruence rings of integers in number fields 
(algebraic extensions of the rationals).
The geometric interpretation is emphasized.
We conclude \S\ref{S:Conclusions} discussing briefly some important topics at hand,
like Frobenius elements and Weil zeros.

\section{Finite Fields: the ``Abstract Way''}\label{S:BG}

We will recall the basic facts about finite fields, as introduced in most standard texts of abstract algebra.
To keep it self-contained, and simple, we use a brief presentation available on the web \cite{Intro-FF}.
See \cite{FF} for additional theoretical details and \cite{Computations-FF} for a computational approach.

\begin{defin}
A finite field is a field which is finite!
\end{defin}
The additive order of the unit $1+...+1=0$ is called the characteristic of the finite field. 
It is always a prime $p$.
For example $F_3=Z/3Z$ has characteristic $p=3$.

Recalling some basic properties are in order.
\begin{prop}
A finite field $F$ of characteristic $p$ has $q=p^n$ elements.
It is a vector space of dimension $n$ over the primary field $F_p$: $n=[F,F:p]$.
\end{prop}

\begin{theorem}\label{T:FF-Existence}
(i) (Existence and Uniqueness) For each $p$ and $n$ there exists a finite field of characteristic $p$
with $p^n$ elements.

(ii) Two such finite fields with the same number of elements are isomorphic.
\end{theorem}
It is therefore natural to denote a {\em generic} finite field as $F_q$, with $q=p^n$, as if it is a specific one.
By abuse of notation, yet well justified by the uniqueness modulo isomorphism,
we write $F_p=Z/pZ$, without further comments (LHS is a ``any'' finite field of characteristic $p$,
while the RHS is the preferred, specific construction of one such field).

The ``standard'' way to construct higher dimensional finite fields with a given number of elements $p^n$,
and of course prescribed characteristic $p$, uses the standard algebraic construction of field extensions
via polynomial rings and their quotients by ideal generated by irreducible polynomials.

We reproduce here the Example 1.88, from \cite{FF}, p.34.
\begin{example}\label{Ex:F3}
Let the prime field be $F_3$.
As an example of the formal process of {\em root adjunction}, 
consider the irreducible polynomial $f(x)=x^2+ x +2\in F_3[x]$. 
Let $\theta$ be a "root" of f; that is, $\theta$ is the residue class $x +(f)$ in $L = F_3[x]/(f )$. 
The other root of $f$ in $L$ is $2\theta+2$, since 
$$f(2\theta+2)=(2\theta+2)^2+(2\theta+2)+2=\theta^2+\theta+2=0.$$ 
We obtain the algebraic extension $L = F_3(\theta)$ consisting of the nine elements
$0,1,2,\theta,\theta+1,\theta+2,2\theta, 2\theta+1,2\theta+2$,
i.e. an instance of $F_{3^2}$.
\end{example}

\section{What {\em are} Number Fields?} \label{S:NumberFields}
The algebraic structure we call {\em field} was first introduced by Dedekind \cite{Dedekind,Stillwell}.
The usual number systems $Q, R$ and $C$ are the traditional examples of fields.
When solving algebraic equations defined by polynomials,
we are ``forced'' to extend our number system, and adjoin formal roots of polynomials as new ``numbers''.
We can treat these either as new {\em symbols}, and {\em construct} the new number system,
for example $C=\{x+iy| x,y\in R, i^2=-1\}$, as real linear combinations of 1 and the symbol $i$
subject to the relation $i^2=-1$, or more formally, in the abstract (algebra), 
as quotients of polynomials modulo the ideal generated by the polynomial defining the relation:
$$C=R[X]/<X^2+1>=\{a+bI | a,b \in R, I=[X]\},$$
Here $[X]$ denotes the congruence class of $X$ modulo the ideal, satisfying the required 
relation: $I^2+1=[X^2+1]=0$ (since $X^2+1\cong 0\ mod \ X^2+1$).

We will call this construction the {\em standard algebraic construction} of a field extension.

Now ``integers'' play a central role in arithmetic, in various rings, and they satisfy the structure of
lattices. 
Initially we may call ``integers'' the subring of a field extension which emerges as a 
corresponding field of fractions,
but field extensions require more care when defining the concept of {\em algebraic integer 
of a field extension}.
\begin{defin}
A {\em latice} $\C{L}$ is a $Z$-submodule of a ring. 
\end{defin}
In particular a lattice is a finitely generated abelian group, 
and can be interpreted and visualized as a ``discrete (finite dimensional) vector space''
(by abuse of language, when there are still relations among generators).

Two good examples of such lattices of algebraic integers are the Gaussian integers and Eisenstein integers.

\subsection{Gaussian Integers}
Complex numbers are a familiar example of field extension of the reals.
To keep the theory algebraic, and to investigate it from an arithmetic point of view, neglect the 
Cauchy reals as non-realistic numbers \cite{Real-fish},
and consider the quadratic extension $Q(\sqrt{-1}$ over the rationals $Q$.
Even better, since these fields are fields of fractions, focus on the extension of integers:
$Z[i]$\footnote{In general the extension of integers might not coincide with the algebraic integers 
of the corresponding field extension.}.

The ring $Z[i]=\{m+in|m,n\in Z\}$ is called the ring of Gaussian integers.
The rational primes $p$ may factor in this larger arithmetic number system:

$\bullet$ $2$ is special, and ``ramifies'' as $2=i^{-1}(1+i)^2$
\footnote{Recall that we have more units here: $1,-1,i,-i$.}.

$\bullet$ $p\cong 1\ mod \ 4$ splits into a product of conjugate primes, like for example $5=(2+i)(2-i)$;

$\bullet$ $p\cong -1\ mod \ 4$ is {\em inert}, i.e. it remains a prime in $Z[i]$; for example $p=3$.

For more facts about Gaussian integers see \cite{Wiki:Gauss}.
For a more technical account, including relations to Galois theory, see \cite{ANT}.

\subsection{Eisenstein Integers}
Similarly, taking a cubic root of unity $\omega$ instead of the 4-th root of unity $i$,
we obtain the Eisenstein integers $Z[\omega]$, with its own primes and classes of rational primes
ramifying ($p=3$), and splitting or being inert, according to a similar condition
$p\cong \pm1$, but this time modulo $3$.
Alternatively, one may look at the analog of Fermat's Two Squares Theorem,
about representing primes $p=m^2+n^2$, except this time we use a different 
quadratic form (norm): $x^2-xy+y^2$, instead of the usual one $x^2+y^2=N(x+iy)$ in $Z[i]$.

For more details, see \cite{Wiki:Eisenstein}.

\subsection{From Number Fields to Finite Fields}
Now the idea for constructing higher dimensional finite fields, is to consider the congruence rings 
of algebraic integers, modulo a prime,
the obvious analog of the construction of primary finite fields $F_p=Z/pZ$.

As a quick example, 
$Z[i]/3Z$ yields $F_{3^2}$, while $Z[i]/(2+i)Z\cong Z/5Z=F_5$.

Besides being a more natural construction, it provides the geometric background for a better understanding
of finite fields as {\em Klein geometries} (Galois fields)
\footnote{... not to mention the connection with Galois Theory,
splitting polynomials and Frobenius elements.}

\section{Lattice Models: the ``Geometric Way''}\label{S:LM}
We will proceed by way of example.
Recall that the primary fields $F_p$ can be constructed as ing quotients $Z/pZ$, where $p$ is a prime 
number, the characteristic.
geometrically, $Z$ can be viewed as a 1-dimensional lattice, or as an infinite oriented graph
\footnote{The interested reader may lookup partial ordered set, POSet for short, too, 
as a generalization.}.

The prime $p$ defines a period, and the covering map $\phi(k)=k\ mode\ p$ is a discrete geometric analog
of the familiar covering map of the circle $x\ mod\ 1$, sometimes used to define angles, sine and cosine.
Algebraically, $\phi$ is a group (ring) homomorphism: the quotient map of the ring $Z$ by the ideal $pZ$.

Now let's consider a 2D-example: the Gaussian integers, as a lattice, modulo a prime ideal $\C{P}$.

Since $Z[i]$ is a principal ideal domain (PID), we need only consider $\C{P}=Z[i]\pi$ with 
Gaussian prime $\pi$ ``sitting'' over a rational prime $p$:
$N(\pi)=p$.

For example $2+i$ is a Gaussian prime over $5$, completely splitting it: $5=(2+i)(2-i)$. 
Recall that other rational primes of the form $p\cong -1\ mod\ 4$ are {\em inert}, i.e. are Gaussian primes too and
$N(p)=p^2$. 

There is also the special case of the {\em ramified prime} $2$ \footnote{Divides the discriminant of the quadratic extension $[Z(i):Z]$.},
which factors with multiplicity: $2=(1+I)^2\cdot (-i)$ 
\cite{Wiki:Gauss}).
\begin{rem}
The factorization may also be written in an initially misleading way as $2=(1+i)(1-i)$,
but $1+i$ and its conjugate $1-i$ are the ``same'' prime, modulo a unit $\pm1, \pm i$.
\end{rem}

Consider the same algebraic quotient map $\phi(z)=z\ mod\ \pi$. 
Since $\pi$ is prime, the quotient ring $K=Z[i]/(\pi)$ is a field of characteristic $p$, i.e. $F_{p^f}$.
The norm $N(\pi)=p^f$ gives the dimension $f=[K:F_p]$.

Excepting the case of the ramified prime $p=2$, we have the following two cases.
For inert (rational) primes $p\cong 3\ mod \ 4$, $\pi=p$ is the only prime over $p$, and $f=2$;
otherwise $p=\pi \bar{\pi}$ splits and $N(\pi)=p$.

\begin{example}
Let $p=5$ and $\pi=2+i$. Then $K=Z[i]/(\pi)$ is a lattice model of $F_5$ (the abstract finite field with $5$ elements).
We can see its canonical residue classes as the Gaussian integers in the fundamental region
of the lattice $\C{L}=\{a\pi+b\bar{\pi}|a,b \in Z\}$, for example with $a,b$ non-negative integers, 
such that $N(z)<p$ (again considering the projection on the integers).
\end{example}

Another example of lattice model, providing an alternative construction to the ``standard''
algebraic extension from Example \ref{Ex:F3}, is the following.
\begin{example}
Consider again $Z[i]$ as a quadratic extension and $p=3$ the rational {\em inert} prime.
Then the quotient lattice $Z[i]/(3)$ has $q=3^2$ elements, representing the finite field $F_{3^2}$.
\end{example}

\section{Applications to Weil Zeros}\label{S:Applications}
There are several topics of Algebraic Number Theory which may benefit from 
the introduction of finite fields as quotients of lattices of algebraic integers:

\quad a) Ramification Theory, in the context of Galois Theory of such extensions;

\quad b) The Frobenius element, as a generator of the Galois group of the corresponding extensions,
controlling the factorization of prime ideals in extensions of number fields;

\quad c) Quadratic Reciprocity using the connection between the Frobenius element and Legendre symbol
in congruence rings of number fields;
finally, 

\quad d) Applications to Weil Conjectures,
and notably to the finite characteristic Riemann Hypothesis via the characteristic polynomial
of a lift of the Frobenius element, having eigenvalues the Weil zeros of 
the Weil polynomial, i.e. the reciprocal of the numerator of the
Hasse-Weil Congruence Zeta Function \cite{Sutherland-EC}.

The first three applications are essentially described in \cite{Ash:ANT}.
In this article we will focus on this later important application to Algebraic Geometry, 
which can be accessed relatively easily, in a computational oriented way,
using for example SAGE as a mathematical software.
In this brief note, we will only point the way. 
For an exposition, see the classical texts, for example \cite{Rosen, Lorenzini};
additional explanations and computations can be found in the lecture notes of the first author \cite{LI:Mat410}.

\subsection{Solving Algebraic Equations over Finite Fields}
Quadratic equations were studied since ancient times, e.g. Appolonius' theory of conic sections.
Replacing the usual number system with finite fields places the problem in the context of Algebraic Geometry.

Following \cite{Lorenzini}, Ch.8, consider the solution $X(F_q)$ of the equation $y^2=x^d+D$ over 
the finite field $F_q$ with $q=p^n$ elements. It is an algebraic {\em affine curve} of {\em degree} $d$.
Denote the corresponding number of elements $N_n$, and the associated {\em congruence zeta function}
\footnote{Conform Weil Conjectures/Deligne Theorem.}
$$Z_{X/Fp}(T)=\frac{P(T)}{1-T)(1-pT)}, \quad P(T)=\prod\limits_{i=1}^{2g} (1-\omega_iT),$$
where $g=(d-1)/2$ is called the {\em genus of the curve}, and $w_i$ are algebraic numbers we
will call the {\em Weil zeros} of the {\em Frobenius polynomial}
$h(u)=u^{2g}P(1/u)$\footnote{$w_i$ are reciprocal of the zeros of the ``Frobenius polynomial'' $P(T)$.}
We will not go in depth explaining the terminology, and just use it to exemplify the relation with 
factorization of primes and lattice models of finite fields.

\begin{example} The cubic ($d=3$) curve $X: y^2=x^3+D$, is an {\em elliptic curve} of genus $g=1$,
which should be pictured topologically (over complex numbers) as a torus 
(when completed with the point at infinity: the {\em projective curve}).
\end{example}
Regarding the fixed prime $p$, whether $F_p$ has $m$-roots of unity or not decides the form of 
$P(T)$ and $N_n$.
In what follows we will assume $m|p-1$, i.e. $F_p$ has $m$-roots of unity
\footnote{Cauchy Theorem for the multiplicative group $(F_q^\times,\cdot)$.}
Then $P(T)=(1-wT)(1-\bar{w}T)$ 
is a quadratic polynomial and the number of affine points is $N_n=p^n-w^n-\bar{w}^n$, where 
$\bar{w}$ denotes complex conjugation \cite{Rosen}, Ch.18\S2, p.302 
(where the $+1$ stands for the point at infinity; see also \cite{Lorenzini}, p.292).
\begin{rem}
Later we will see how Weil zeros $w_i$ are related to Gauss and Jacobi sums,
which are valued in the cyclotomic numbers of roots of unity of order
$l(l-1)$ and $l-1$ respectively, if $l=m|p-1$ is a prime.
\end{rem}
Part of {\em Weil Conjectures} \cite{WC} is that $w\bar{w}=p$, 
i.e. the {\em Riemann Hypothesis} holds in finite characteristic \cite{Lorenzini}.
Moreover, introducing the {\em defect} $a_p=w+\bar{w}$, 
$P(T)=1-a_pT+pT^2$.
\begin{rem}
The coefficients of the Betti polynomial $P(T)$ are related to Weil zeros as a consequence 
of a deeper connection with the {\em characteristic polynomial of the Frobenius element} $h(u)$:
$a_p=Tr(Frob)$, $p=det(Frob)$.
\end{rem}
\begin{example}
The elliptic curve $m=3$ \cite{Rosen}, p.306, has $N_1=p+w+\bar{w}$, 
where the Weil zeros split the prime $p=w\bar{w}$ in the cyclotomic extension $Z(\zeta_3)$
of Eisenstein integers (assuming $3|p-1$).
In terms of primary primes $\pi,\bar{\pi}$, $\pi\cong 2 \ mod \ 3$ (associated to $w, \bar{w}$),
we have \cite{Rosen}, Th.4, p.305 (affine points; $6=l\cdot (l-1)$ with $l=3$):
$$N_p=p+2 Re(\rho(4))\pi, \quad \rho(x)=\Big(\frac{4D}{\pi}\Big)_6\pi.$$
As a concrete example take $D=1$.

If $p=13$ then $\pi=-1+3\omega$ is a primitive prime, and together with $\bar{\pi}=-1+3\omega^2$ split $p$:
$$p: \ 13=(-1+3\omega)(-1+3\omega^2)\ : \pi\bar{\pi}.$$
Since $\rho(4)=\omega^2$, the Weil zero is $w=-\omega\pi$, associated to $\pi$ 
(Units: $\{\pm1, \pm\omega, \pm\omega^2\}$).

Now the number of affine (finite) points in $F_p$ is:
$$N_1=13+2 Re(\omega \pi)=13+2(\omega^2+\omega)=13-2=11,$$
consistent with a direct check and counting argument:
$$X(F_{13})=\{(4,0), (10,0), (12,0), (0,\pm1), (2, \pm3), (5,\pm3), (6,\pm3)\}.$$
\end{example}
\begin{example}
As another example consider $D=5$ and $p=19$.
Note that $\pi=5+3\omega$ is primary and splits $p$:
$$19=(5+3\omega)(5+3\omega^2), \ 5+3\omega\cong 2\ mod\ 3.$$
From the above formula we obtain the number of points (\cite{Wang},p.8)
\footnote{We need Jacobi sums for this: see \S \ref{S:JS}.}:
$$N_1=p+\pi+\bar{\pi}=19+5+\omega+5+\omega^2=26.$$
Then $a_p=-2 Re(\pi)=-7$ and Weil zeros are $w=-\pi$ and its conjugate:
$$P(T)=(1-wT)(1-\bar{w}T)=1+7T+19T^2, \quad N_1^{proj}=P(1)=27.$$
\end{example}
\begin{rem}
One may use {\em SAGE} (recently renamed as {\em CoCalc}) \cite{LI:SAGE-ANT} 
to conveniently compute {\em Dirichlet characters} (like $\rho$ above), and Jacobi sums, 
which are instrumental in computing the number of points.
\end{rem}
Well, what does this problem of counting the number of solutions,
with its associated congruence zeta function,
have to do with lattice models of finite fields!?
For this we need to recall some facts about the Galois group of an extension, 
and the relation with the {\em Frobenius element},
which will turn out to be present as the numerator of the zeta function.

\subsection{Frobenius Element}
Following \cite{Yott}, consider a number field extension $Q(\xi)/Q$ whic is Galois,
and how a rational prime $p$ decomposes in it, 
with $\pi$ such a prime factor (assuming a principal ideal domain case for simplicity).
Then the Galois group of the number field extension $Gal(Q(\xi)/Q)$ 
is related to the Galois extension of the corresponding (lattice models of) finite fields:
$$1\to I(\pi\to D\to Gal(F_q/F_p)\to 1,$$
where $D(\pi/p)$, the {\em decomposition group} consists in Galois automorphisms
preserving the ideal generated by $\pi$, each of its elements therefore inducing
an automorphism of the corresponding finite fields extension $Gal(F_q/F_p)$,
in a surjective manner.

If $p$ is unramified, then the kernel (the {\em inertia group}) is trivial,
and the above surjection becomes an isomorphism.
Then one can ``pull-back'' the Frobenius automorphism $x\mapsto x^p$
of $Gal(F^q/F^p)$,
where we recall that $F^q=Z[\xi]/\pi$ and $F_p=Z/p$ are lattice models of finite fields
constructed in number fields {\em viewed (embedded) as subfields of the complex numbers}.

\begin{defin}
In the context above ($\pi$ prime in $Z[\xi]$ over the unramified rational prime $p$), 
the {\em Frobenius element} $Frob_p^\pi\in Aut_Z(Z[\xi])$ is the unique
Galois automorphism which induces the Frobenius automorphism $Fr(x)=x^p$
in the finite field extension $F_q/F_p$, of lattice models.
\end{defin}

At this stage the Frobenius elements may depend on the choice of prime $\pi$ over $p$.
But these Frobenius elements are conjugate to each other,
so if the Galois group is Abelian, then the Frobenius element is unique,
and will be denoted by $Frob_p$.

\begin{example} \label{Ex1}
Consider $Q[i]$, with $i$ a forth root of unity, and its Gaussian integers $Z[i]$.
The only ramified prime is $2$; otherwise $p\cong 1$, $mod \ 4$ or course, splits,
or $p \cong -1$ is inert.

The decomposition group $D(p)$ is trivial in the split and ramified cases,
and equals $G=Gal(Q(i):Q)\cong Z_2$ (multiplicative group $\{-1,1\}$) otherwise.

Thus the Frobenius element is $1$ when $p\equiv 1$ and $-1$ otherwise, 
i.e. $Frob_p=(\frac{-1}{p})$ is given by the {\em Legendre symbol}
(the unique multiplicative character of order $2$).

Alternatively, we can compute the {\em lift to $Z[i]$ of the Frobenius} $x^p$, 
from the ``abstract'' setup, using our lattice model,
to the ring of algebraic integers\footnote{It is enough to consider the extension of $Z$,
and not the full algebraic closure in $Q(i)$, which incidenltally, here coincide.}:
$$ (a+ib)^p\equiv a+b i^p \quad mod \ pZ[i].$$
Since $i^p=(\frac{-1}{p})i$ ``on the nose'', i.e. not just $mod \ p$, we conclude that
the lift has the following closed formula in terms of the 
multiplicative quadratic residue character $\rho_2(z)=(z/p)$: 
$$Frob_p(a+ib)=a+\left(\frac{-1}{p}\right)b,
\quad Frob_p=\sigma^k, \ k=\left(\frac{-1}{p}\right)=(-1)^{\frac{p-1}{2}}.$$
The last form was written in terms of the generator of $G$, here complex conjugation.
\end{example}

\begin{example}\label{Ex2}
The above example can be generalized to quadratic extensions $Q(\sqrt{d})$,
where $d$ is square free (\cite{Yott}, p.3).
The Frobenius element, in $Z/pZ^\times$, is  
$Frob_p=(d/p)$, so that $Frob_p(a+\sqrt{d} b)=a+(-1/p)\sqrt{d}b)$,
i.e. $Frob_p=\sigma^{ord(-1/p)}$ as before.
\end{example}

\begin{example}\label{Ex3}
In the cyclotomic case $Q(\xi_n)$, the primes that ramify are those which divide $n$.
The Galois group is isomorphic to the multiplicative group of roots of unity,
and therefore isomorphic to $(Z/nZ^\times, \cdot)$, 
with a Galois element $\sigma_m:\xi\mapsto \xi^m$, with $m\in Z/nZ^\times$ relatively prime to $n$.
 
In the non-ramified case $p\in Z/nZ^\times$, the Frobenius element is, 
again as expected: 
$$Frob_p(\sum_{k=0..n-2} c_k\xi_n^k)=\sum c_k Frob_p(\xi_n)^k, \quad Frob_p(\xi)=\xi_n^p.$$
As another quadratic extension example, consider $Q(\omega)$, corresponding to a cubic root of unity $\omega^3=1$, 
and its Eisenstien integers $Z[\omega]$.
Then the corresponding Frobenius element is, similarly to the Gaussian integers case:
$$Frob_p(a+\omega b)=a+\rho_3(p)\omega b.$$
\end{example}
\begin{rem}
At this stage, one may further look into the correspondence 
between how the prime $p$ factors into $Z[\xi]$, and how 
the primitive polynomial $f(x)$, of $\xi$  factors in $F_p[x]$,
reflecting the commutativity of the diagram \ref{Diagram1} from the introduction.
\end{rem}

It is conceptually important to piece together these Frobenius elements as a map
depending on the prime $p$, called the {\em Artin map}: $Frob: p\mapsto Frob_p$
\footnote{In a more general setup \cite{Ash:ANT}, Ch.5, 
$Frob(p)=(\frac{L/K}{P})$ is called the {\em Artin symbol}.}.
For cyclotomic extensions, If we identify the Galois group $Gal(Q(\xi_n)/Q)$ with
$(Z_n^\times,\cdot)$, then the Artin map is simply the ``identity'' map:
$$p\mapsto Frob_p=p \ mod \ n, p\in Z_n^times.$$
For example, with $m=4$ and $p$ an odd prime,
the Galois group is generated by conjugation $G=<\sigma>$,
since $Z_4^\times\cong Z_2$,
and $Frob(p)=p \ mod \ 4$ as an element of $Z_4^\times$.
This is essentially the Legendre symbol $\frac{p}{4}$,
when identifying $Z_4$ with the 4-th roots of unity, via exponentiation
(the Galois group identification).

Once we know the Frobenius element,
its characteristic polynomial can be computed easily:
$$P(Forb_p)(t)=det(Frob_p-t I).$$
For example, in the cyclotomic setup, with $m=4$ (Gaussian integers),
the matrix of $Frob_p=\sigma^{(-1/p)}$ in the basis $1, i$ is: 
$$Split: \ p\equiv 1\ mod \ 4: \quad Frob_p=
\begin{bmatrix} 1 & 0 \\ 0 & -1 \end{bmatrix}$$
$$Inert: \ p\equiv-1 \ mod \ 4: \quad Frob_p=
\begin{bmatrix} 1 & 0 \\ 0 & 1 \end{bmatrix}
$$
and the characteristic polynomials are, respectively:
$$ P(T)=(1-T)^2, \quad P(T)=1-T^2.$$
Similarly, for a quadratic extension like for $m=3$ (Eisenstein integers), 
the matrices of the Frobenius elements $I$ and $\sigma$, 
and their matrices are essentially the same 
(but computed in a different basis $1,\omega$).

Now let's see how the Frobenius element, 
or rather its lift and the corresponding characteristic polynomial
is related to the Hasse-Weil congruence zeta function.

\subsection{Weil Zeros and Jacobi Sums}\label{S:JS}
We will only document the facts with an example, following \cite{Silverman,Ang,Rosen, Lorenzini},
and leave the general case for a separate study.

Let $y^2=x(x^2+1)$ define an elliptic curve over $F_q$.
Since the RHS of its defining equation $f(x)$, splits in $Z[i]$, we will work with Gaussian integers
in the number fields side of the ``picture''.

For $p\cong 3\ mod \ 4$, the prime is inert in $Z[i]$,
which corresponds to the factor $x^2+1$ being irreducible in $F_p$ 
and the Frobenius element complex conjugation.

Theorem 5, \cite{Rosen}, p.307, with $D=-1$, yields the number of {\em projective} points
\footnote{The $+1$ stands for the point at infinity.},
according to the type of prime:
$$Inert:\ N_p=1+p, \quad Split: \ N_p=p+1- 2 Re(\rho_4(-1)\pi,$$
where $\pi$ is a primary prime splitting $p$ and $\rho_4$ is a character of order $4$.

We will focus on the split case $p\cong \ 1\ mod\ 4$ 
(Ramification Theory parameters: $g$=2, $e=1$, $f=1$).

To have a ``nice'' description of the lift of Frobenius $Fr_p$ on $C$ preserving our curve, 
and not some ``deformation'' of identity (the Frobenius element)
$(x,y)\mapsto (x^q+f(x,y), y^q+g(x,y))$ \cite{Ang}, p.10,
we use Weierstrass coordinates.
The elliptic curve is then the quotient of $\B{C}$ by our lattice $\Lambda=Z[i]$ of (Gaussian) algebraic integers:
$$e(\tau z): (\B{C},+)\to (\B{C}^\times,\cdot), \quad \B{C}/Z[i]\cong E(\B{C}),$$ 
where here $\tau=2\pi i$.
Then the Frobenius lift $Fr_p(z)=z\cdot c$ 
is multiplication by some lattice element $c=a+ib\in Z[i]$ \cite{Ang,Silverman}.
\begin{rem}
Alternativelly, we could lift the Frobenius to the p-adic completion, and 
taking advantage of Hasse principle for finding the above ``perturbations'' 
$f(x,y)$ and $g(x,y)$.
\end{rem}

If the curve is defined by a polynomial in the powers of the variables (Weil curves), 
e.g. Riemann surfaces $y^2=x^s+D$ (\cite{Lorenzini} p.292) and 
Fermat curves $x^m+y^m=z^m$ (\cite{Rosen,Lemmermeyer}) , 
then Jacobi sums provide a powerful tool to compute the number of points.

Then $N_p=1+p-a_p$, with the {\em defect} given by the Jacobi sum $a_p=2 Re J(c_2,c_4)$,
which also yields the Weill zeros $w,\bar{w}$ of the (reciprocal of) ``Betti polynomial''
\footnote{The numerator of the Zeta function is a local L-function
having a cohomological interpretation.}:
$$Z_p(T)=\frac{L_p(T)}{(1-T)(1-pT)},\quad L_p(T)=(1-wT)(1-\bar{w}T), \quad w\bar{w}=q.$$
Then $w=-J(c_2,c_4)$ is primary \cite{Rosen}
and our lift of Frobenius is given by $c=w=\pi$,
conform with \cite{Silverman},
with $a_p=Tr(Fr_p)$ and $p=det(Fr_p)$ 
(Riemann Hypothesis, part of the Weil Conjectures; 
see also \cite{Southerland-EC}, Lecture \#8, Hasse's Theorem):
$$ CharPoly(Fr_p): \quad det[Fr_p-u\ Id]=u^2-a_pu+p, \ u=1/T.$$
Rewriting the number of points as in \cite{Ang} $N=q-2d\sqrt{q}+1$,
and the Weil zero as $w=e^{i\theta}\sqrt{q}$, 
one may interpret the ``Betti coefficient'' $d=Re(e^{i\theta}$
via a comparison with the Jacobi sum, as phase of the 
2-cocycle of the Fourier coefficients of the 
Dirichlet characters (Gauss sums) ... but this is another story!

\begin{rem}
A similar discussion applies to our previous example of elliptic curve 
$y^2=x^3-D$ (\cite{Rosen}, p.304; \cite{Wang}, p.7), 
with Eisenstein integers replacing Gaussian integers.
\end{rem}
\begin{rem}
For higher dimensional extensions $F_{p^f}$, needed when the genus of the curve 
excedes $g=1$, can be implemented via cyclotomic extensions $Z[\xi_m$,
such that the dimension $n=\phi(m)$ factors as $g\cdot f$ with
$f=ord(p)$ is the order of multiplication by $p$ in $Z/mZ^\times$
and $g$ the ramification genus of the prime $p$.
\end{rem}

\section{Conclusions and Further Developments}\label{S:Conclusions}
There are various styles of teaching (and designing) Abstract Algebra.
We attempted to plead that, in the case of finite fields,
the abstract approach to the introduction of the algebraic structure (``axiomatic''/top-down design),
can be supplemented by the specific construction we call {\em lattice models}, 
which introduces the number fields first, as more ``familiar'' to the student used to solve polynomial equations, 
and presenting $F_{p^n}$ as a congruence ring, in perfect analogy to the way we introduce 
the primary finite fields $F_p=Z/pZ$.

The bonus is some extra intuition, but more importantly a rich {\em geometric framework} 
for bridging and interpreting other abstract algebra concepts, like Galois Groups, Frobenius elements,
paving the road towards understanding {\em General Reciprocity Laws} the ``right way'' \cite{Ash_and_Gross}.

As another more advanced application, Weil Conjectures can be understood not
in their natural ``habita't'' of abstract Algebraic Geometry,
but in the more geometric and topological context of complex manifolds,
by using lattice models of finite fields.
Then the Frobenius element of the number fields incarnation of the finite fields is their Frobenius element,
without the need of a lift. 
Then the numerator of the Weil Congruence zeta function
is the characteristic polynomial of the lift of the Frobenius element, allowing to count numbers of solutions
without the use of a Weil cohomology (e.g. Grothendieck's approach via l-adic cohomology).

Yet if one wishes, Lefshetz formula, as well as algebraic topology/geometry technics,
such as Riemann-Roch/Hurwitz Th.,  may be used on this characteristic zero ``side''
of number fields and lattice models of finite fields, 
in the natural and familiar framework of the complex numbers.

On the concrete side, computer programs are included for computer explorations of the 
presented topics of Algebraic Number Theory \cite{LI:SAGE-ANT}.

\appendix

\section{SAGE / CoCalc Programs}
Programs for computing the Weil zeros of $EC: y^2=x^3+D$ are available 
from \cite{LI:SAGE-ANT}.
They can be easily adapted to other cases, for example to Riemann Surfaces or Fermat Curves. 
The programs can also be used to compute Jacobi sums, and for other Algebraic Number Theory studies
using technology.



\end{document}